\documentclass[12pt]{article}  

\usepackage{amsmath,amsxtra,amssymb,latexsym, amscd}
\usepackage[mathscr]{eucal}

\begin{document}

\title{ \huge\bf Algebroids and Jacobian conjecture 
 }

\def\1{\rule{0cm}{0cm}} \def\qd{\rule{3mm}{3mm}} \def\BB{$\bullet$}
\def\bz{\boldsymbol{z}} 
\def\bx{\boldsymbol{x}} \def\by{\boldsymbol{y}} \def\bS{\boldsymbol{S}}
\def\I{{\cal I}} \def\A{{\cal A}} \def\C{\mathbb{C}}
\def\bv{\boldsymbol{v}}
\def\bn{\boldsymbol{n}} \def\by{\boldsymbol{y}} \def\bS{\boldsymbol{S}}
\def\H{{\mathcal H}} \def\U{{\cal U}} \def\D{{\cal D}}
\def\bc{\boldsymbol{c}} \def\bZh{\boldsymbol{Z}}
\def\ba{\boldsymbol{a}}\def\bu{\boldsymbol{u}}
\def\br{\boldsymbol{r}}\def\bzero{\boldsymbol{0}}

\def\R{\mathbb{R}} \def\bw{\boldsymbol{w}}
\def\C{\mathbb{C}}

\author{
{\bf Dang Vu Giang}\\
Hanoi Institute of Mathematics\\
Vietnam Academy of Science and Technology\\
18 Hoang Quoc Viet, 10307 Hanoi, Vietnam\\
       e-mail: dangvugiang@yahoo.com\\
\1\\
}
\maketitle 

\noindent {\bf Abstract.} We prove the  injectivity of any kellerian mapping using the fact that the holomorphy domain of an algebroid defined via an irreducible polynomial may contain singular points of that polynomial. The famous Jacobian conjecture is true. 

\bigskip



\section{Introduction}

\bigskip
\par\noindent
Kellerian mapping is a polynomial mapping $\Phi:\C^n\to\C^n$   whose  Jacobian is 1. The famous Jacobian conjecture  \cite{essen}  says that  any kellerian mapping   is a polynomial automorphism. Thousands of authors are interested in proving this conjecture but unsuccessfully.  In the virtue of Galois theory and Nullstellensatz, any injective polynomial mapping is a polynomial  automorphism. Therefore,  it is enough to prove that any kellerian mapping is injective \cite{essen}. My idea is based on the holomorphy domain of  algebroids, which may be holomorphic in any small neighborhood of singular points. For example, the algebroids $T(z)^2-z^2-z^3=0$ are holomorphic  at $z=0$, (holomorphic in a small neighborhood of this point).  On the other hand, my arguments do not work in famous Pinchuk real map, because the Jacobian of Pinchuk map may be becoming 0 at complex coordinates. In fact, the Pinchuk map is not open in $\C^n$. Moreover, my arguments work very well for any field $\C_p$ of complex $p-$adic numbers, because we do not use the local compactness of $\mathbb{C}$  the field of complex numbers. It is well known that, if the Jacobian conjecture is true for $\C$ then it should true in any $\C_p$. Recall that algebroid is a continuous function of several complex variables, which satisfies a polynomial equation. For example, does  exist a holomorphic function $T$ around the origin which satisfies the equation $T(z)^2-z^2-z^3=0$. The two variables polynomial $t^2-z^2-z^3$ is irreducible.

\section{Main Result}

\bigskip
\noindent
{\bf Theorem. } {\it  Any locally injective polynomial   mapping $\Phi:\C^n\to\C^n$ 
  is injective. Therefore, the famous Jacobian conjecture is true.}
  
\bigskip
\noindent
{\bf Proof. }   
Clearly, the field $\C(\boldsymbol{z})$ is a finite extension of the field $\C(\Phi)$. Hence, any $z_j$ has a minimal polynomial over $\C(\Phi)$. For example, let 
$f(\Phi,z_1)={{\alpha }_{0}}\left( \Phi \right)z_{1}^{m}+{{\alpha }_{1}} \left( \Phi \right) z_{1}^{m-1}+\cdots +{{\alpha }_{m}}\left( \Phi \right)$ denote  the minimal polynomial of $z_1$. Clearly, this polynomial is irreducible in $\C(\Phi)[z_1]$.  Moreover, 
$${{\alpha }_{0}}(\Phi\left( \boldsymbol{z} \right))z_{1}^{m}+{{\alpha }_{1}}(\Phi\left( \boldsymbol{z} \right))z_{1}^{m-1}+\cdots +{{\alpha }_{m}}(\Phi\left( \boldsymbol{z} \right))=0,\qquad \bz\in\C^n.
\eqno(1)$$ 
 Consider the differential field  $\mathcal C$ generated by $\C(\bz)=\C(z_1,z_2,\cdots,z_{n})$
and symbols $\bZh=(Z_1,Z_2,\cdots,Z_{n})$ with the reduced rule $\Phi(\bZh)={\boldsymbol z}$. Clearly, the extension ${\mathcal C }/\C(\bz)$ is finite. 
Moreover,
${{\alpha }_{0}}\left( \boldsymbol{z} \right)Z_{1}^{m}+{{\alpha }_{1}} \left( \boldsymbol{z} \right) Z_{1}^{m-1}+\cdots +{{\alpha }_{m}}\left( \boldsymbol{z} \right)$ is  the minimal polynomial of $Z_1$. 
Then 
$${{\alpha }_{0}}\left( \boldsymbol{z} \right)Z_{1}^{m}+{{\alpha }_{1}} \left( \boldsymbol{z} \right) Z_{1}^{m-1}+\cdots +{{\alpha }_{m}}\left( \boldsymbol{z} \right)=0.\eqno(2)$$  
On the other hand, $\Phi$ is an open mapping of $\C^n$. 
This means that if the equation system $\Phi({\boldsymbol u})=\boldsymbol{a}$ has a root then 
 in a small neighborhood of $\boldsymbol a$, there is a holomorphic mapping $\bZh=(Z_1,Z_2,\cdots,Z_{n})$ such that $\Phi(\bZh(\boldsymbol{z}))=\boldsymbol{z}$ and $\bZh(\boldsymbol{a})={\boldsymbol u}$. (This argument does not work for Pinchuk map, because the Pinchuk map is not open in $\C^n$.)  Consequently,
 $$
{{\alpha }_{0}}\left( \boldsymbol{z} \right)Z_{1}(\boldsymbol{z})^{m}+{{\alpha }_{1}} \left( \boldsymbol{z} \right) Z_{1}(\boldsymbol{z})^{m-1}+\cdots +{{\alpha }_{m}}\left( \boldsymbol{z} \right)=0,
\eqno(3)$$  for  ${\boldsymbol z}$ in a small neighborhood of ${\boldsymbol a}=\Phi({\boldsymbol u}).$
Moreover, $Z_{1}$ is holomorphic in a small neighborhood of  any point $\bz\in\Phi(\C^n).$ 
Let $$f(\bz,t)=\alpha_0(\bz)t^m+\alpha_{1}(\bz)t^{m-1}+\cdots+\alpha_m(\bz).$$
Then $f(\bz,t)$ is an irreducible polynomial in $\C(\bz)[t]$. 
Now we decompose 
\[f(\bz,t)=\alpha_0(\bz)(t-T_1(\bz) )(t-T_2(\bz) )\cdots(t-T_m(\bz))
\] in an algebraic closure of $\C(\bz)$. Clearly, $T_j$ all are distinct  holomorphic functions around $\ba\in\Phi(\C^n)$ and satisfying $f(\bz,T_j(\bz))=0$.
Now let $U$ denote the zero set of 
\[\frac{\partial}{\partial t}f(\Phi(\bz),z_1).
\] This set is closed and may be empty. If $U$ is non-empty then it is a closed hyper-surface in 
$\C^n$  $(n>1)$. Moreover,
$$\nabla \alpha_0( \Phi (\bz))z_1^m+\nabla \alpha_{1}( \Phi (\bz))z_1^{m-1}+\cdots+\nabla \alpha_m( \Phi (\bz))=\bzero$$ for all $\bz\in U$.
Let $V$ denote the zero set of $\alpha_0( \Phi (\bz))$.
We have at least two $T_j$'s are equal on $\Phi(U\setminus V)$, say $T_1\circ \Phi(\bz)=T_2\circ \Phi(\bz)=z_1$ for all $\bz\in U\setminus V.$   Let $n>1$ and take gradients of $T_j$ for $j=1,2$. We have
$\nabla T_1\circ \Phi(\bz)\Phi'(\bz)=\nabla  T_2\circ \Phi(\bz)\Phi'(\bz)=(1,0,\cdots,0)$ for all $\bz\in U\setminus V$. Here, $\Phi'(\bz)$ denotes the Jacobian matrix of $\Phi$ at $\bz$, which is invertible. Hence, 
$\nabla T_1\circ \Phi(\bz)=\nabla T_2\circ \Phi(\bz)$ for all $\bz\in U\setminus V$.  
Similarly, all derivatives of $T_1$ and $T_2$ are the same at any point $\ba=\Phi(\bz)\in\Phi(U\setminus V)$.   Consequently, the Taylor series of $T_1$ and $T_2$ at any $\ba\in\Phi(U\setminus V)$  are the same. This contradicts the fact that $T_1$ and $T_2$ are two distinct holomorphic functions. 
Hence,  for  all  $\bz\in\C^n$ satisfying 
\[m{{\alpha }_{0}}\left( \Phi({\boldsymbol z}) \right)z_{1}^{m-1}+(m-1){\alpha }_{1} \left( \Phi({\boldsymbol z}) \right) z_{1}^{m-2}+\cdots +{{\alpha }_{m-1}}\left(\Phi({\boldsymbol z}) \right)=0. \eqno(4)\]  we should have
$\alpha_0( \Phi (\bz))=0$,
 The  polynomial in (4) has degree in $z_1$ higher than the degree of  $\alpha_0( \Phi (\bz))$  so we have
$$m{{\alpha }_{0}}\left( \Phi({\boldsymbol z}) \right)z_{1}^{m-1}+(m-1){\alpha }_{1} \left( \Phi({\boldsymbol z}) \right) z_{1}^{m-2}+\cdots +{{\alpha }_{m-1}}\left(\Phi({\boldsymbol z}) \right)=c\alpha_0( \Phi (\bz))^\ell
\eqno(4^*)
$$  ($c$ denotes a non-zero constant).
Now put  $A_j=\alpha_j\left(\Phi({\boldsymbol z}) \right)$ we have two algebraic equations for $z_1$ from (1) and ($4^*$)
$$A_{0}z_{1}^{m}+A_{1} z_{1}^{m-1}+\cdots +A_{m}=0.\eqno(5)$$
$$mA_{0}z_{1}^{m-1}+(m-1)A_{1} z_{1}^{m-2}+\cdots +A_{m-1}=cA_0^\ell\eqno(6)$$
If $m>1$, equation (6) has at most $m-1$ roots $z_1\in\C$,  but by (6) equation (5) has 
exactly $m$ roots in $\C$ if $A_0\ne0$. Hence, $A_0=0$ and we have $\alpha_0\left( \Phi({\boldsymbol z}) \right)=0$ for all $\bz\in\C^n$. Taking gradient we have
$\nabla \alpha_0\left( \Phi({\boldsymbol z}) \right)\Phi'(\bz)=0$. Here, $\Phi'(\bz)$ denotes the Jacobian matrix of $\Phi$ at $\bz$, which is invertible. Hence, $\nabla \alpha_0\left( \Phi({\boldsymbol z}) \right)=0$. Similarly, all derivatives of $\alpha_0$ at $\ba=\Phi({\boldsymbol z})\in\Phi(\C^n)$ are zero. Therefore, $\alpha_0=0$ and this contradicts the fact that the degree of $z_1$ over $\C(\Phi)$ is $m$.
 Therefore, $m=1$  and  $\Phi$ is injective. The proof is now complete.

\bigskip
\par\noindent
{\bf Acknowledgement.} 
Deepest appreciation is extended towards 
 the NAFOSTED  (the National Foundation for Science and Techology Development in Vietnam) for the financial support.

\end{document}